\documentclass[12pt]{article}
\usepackage{bezier,latexsym}
\usepackage{epsfig,afterpage}
\usepackage{amsmath, amssymb}
\usepackage{psfrag}
\usepackage{color}
\catcode`\@=11 
\@addtoreset{equation}{section} 
\catcode`\@=12 
\newcommand\real{{\rm I\! R}}

\newtheorem{theorem}{Theorem}[section]

\newcommand\grad{\mathop{\rm grad}\nolimits}
\newcommand\divv{\mathop{\rm div}\nolimits}

\newcommand\ar{\!\!\!\!}

\def\bql#1{\begin{equation}\label{#1}}
\def\bq{\begin{equation}}
\def\eq{\end{equation}}
\def\eref#1{{\rm (\ref{#1})}}
\usepackage[colorlinks=true, linkcolor=black, citecolor =black, pagebackref=false, bookmarks=false, pdfpagemode=UseNone]{hyperref}

\usepackage{array,multirow}
\usepackage{hhline}
{
\usepackage[table]{xcolor}

\begin{document}
\title{On a boundary integral solution of a lateral planar Cauchy problem in elastodynamics}

\author{Roman Chapko\thanks{Faculty of Applied Mathematics and Informatics, Ivan Franko National University of Lviv, 79000
Lviv, Ukraine},
B. Tomas Johansson\thanks{Mathematics, EAS, Aston University, B4 7ET Birmingham, UK}  and
Leonidas Mindrinos\thanks{Johann Radon Institute for Computational and Applied Mathematics (RICAM), Austrian Academy of Sciences, Linz, Austria}
}

 \date{}
 \maketitle

 \begin{abstract}
A boundary integral based method for the stable reconstruction of missing boundary data is presented for the governing hyperbolic equation of elastodynamics in annular planar domains. Cauchy data in the form of the solution and traction is reconstructed on the inner boundary curve from the similar data given on the outer boundary. The ill-posed data reconstruction problem is reformulated as a sequence of boundary integral equations using the Laguerre transform with respect to time and employing a single-layer approach for the stationary problem. Singularities of the involved kernels in the integrals are analysed and made explicit, and standard quadrature rules are used for discretisation. Tikhonov regularization is employed for the stable solution of the obtained linear system.  Numerical results are included showing that the outlined approach can be turned into a practical working method for finding the missing data.  
\end{abstract}

{\small {\it Keywords:} Boundary integral equations; Cauchy problem; Elastic equation; Elastodynamics; Laguerre transformation; Nystr\"om method; Single-layer potentials;  Tikhonov regularization.}

\section{Introduction}
We assume that we have a two-dimensional physical body modelled as a doubly connected domain $D$ in $\real^2$. This domain has two simple closed smooth boundary curves $\Gamma_1$ and $\Gamma_2$, with $\Gamma_1$ lying in the interior of $\Gamma_2$. Consider the Cauchy problem
 in $D$ for the hyperbolic elastic equation, that is 
\bql{h_eq}
\frac{{\partial^2 u}}{{\partial \,t^2}} = \Delta^* u ,
\quad \text{\rm in }
 D \times ({0,\infty} ),
\eq
with the Lam\'e operator defined by $\Delta^* =c_s^2\Delta+ (c_p^2-c_s^2) \grad  \mbox{div}$, supplied with the homogeneous initial conditions 
\bql{h_ic}
\frac{\partial u}{\partial t}(\cdot,0) = u(\cdot,0) = 0, \quad \text{\rm in } D
\eq
and the boundary conditions
\bql{h_bc}
u= f_2 ,\quad  \text{\rm on } \Gamma_2 \times (0,\infty),\quad Tu=g_2 ,\quad \text{\rm on } \Gamma_2 \times (0,\infty),
\eq
where $f_2$ and $g_2$ are given and sufficiently smooth functions and $T$ the traction operator 
\bql{stress}
Tv =\lambda\divv v\, \nu + 2\mu\,(\nu\cdot \grad)\,v + \mu \divv
(Qv) \,Q\nu ,
\eq  
with the outward unit normal vector $\nu$ to the boundary of $D$, and the unitary matrix 
$$
Q=
\begin{pmatrix}
\phantom{-}0 & 1\\
-1 & 0
\end{pmatrix}.
$$
Here, the velocities $c_s$ and $c_p$ have the following form 
$$
c_s=\sqrt{\frac{\mu}{\rho}}, \qquad c_p=\sqrt{\frac{\lambda+2\mu}{\rho}},
$$
where $\rho$ is the density, and $\lambda$ and $\mu$ are the Lam\'e constants. The components of a generic point $x$ of $\real^2$ is written as a column vector, and we use the notation $x^\top=(x_1,x_2)$ with $x^\top$ meaning the transpose.

The governing equation~(\ref{h_eq}) occurs in elastodynamics and models planar elastic waves, with $u$ being the displacement (small deformations only), in an isotropic elastic medium, and $c_p$ and $c_s$ are the velocities of what is known as the pressure and shear waves, respectively (sometimes also called the primary and secondary waves, note that $c_p>c_s$ for positive Lam\'e constants). For a derivation of the governing equation in the field of elasticity, see for example~\cite[Chapt. III, \S 22]{LandauLifshitz}. Applications of the present work can thus be found in structural engineering to forecast from incomplete data vibrations (typically standing waves) in beams supporting buildings and bridges, and also in seismology. Principal materials where disturbances or vibrations occur and where the governing equation can be applied include solid and dense liquid media, like metals, rocks and water.  

Boundary data of the form~(\ref{h_bc}) is not specified on all of the boundary and this typically leads to an ill-posed problem termed a Cauchy problem. Existence of a solution to the boundary problem, where displacement is imposed on one boundary part and the traction on the other is a classical well-posed problem, existence and uniqueness of a weak solution follows for example from energy integral methods or from results on abstract (non-linear) hyperbolic equations such as~\cite{Hughes, Dafermos} covering also non-linear models (for more direct results focusing on the linear case, see~\cite{Gurtin,Knops}). From a solution to the well-posed problem one can generate data to the ill-posed one and see that it will be solvable for some classes of given functions. Uniqueness of a solution to the Cauchy problem follows from~\cite[Chapter~5]{Tataru}. 

Although existence of a solution to~(\ref{h_eq})--(\ref{h_bc}) can be assumed the solution will in general not depend continuously on the data. Such problems and methods for reconstructing the missing data in a stable way are well-studied for elliptic and parabolic equations, see~\cite[Chapt.~3]{Isakov} for an introduction. However, for hyperbolic equations ill-posed lateral Cauchy problems are less studied, for some works, see~\cite{Yamamoto, ChaJoh, Klibanov} and also~\cite{Hasanov}. An overview of inverse problems in elasticity is given in~\cite{Bonnet}. For some works on inverse problems for the related system of thermoelasticity, see~\cite{Bellassoued, Isanov, Karageorghis, Kozlov09, van_Bockstal}.

In~\cite{ChaJoh} a boundary integral equation based method is presented for both parabolic and hyperbolic ill-posed lateral Cauchy problems. That work follows a string of works by the same authors, where the solution to Cauchy problems are represented in terms of boundary integrals, see further~\cite{ChaJoh16}. Moreover, in~\cite{ChMi} a boundary integral equation method is given for the direct Dirichlet problem for the time-dependent elastic equation in an unbounded two-dimensional domain. In both works~\cite{ChaJoh,ChMi} the Laguerre transformation is employed to reduce the problem to a sequence of stationary problems. We shall combine the ideas of~\cite{ChaJoh,ChMi} to present a stable boundary integral based method to reconstruct data on $\Gamma_1$ given (\ref{h_eq})--(\ref{h_bc}). This involves an investigation of the singularities of the kernels in the boundary integrals and a further technical part to derive a suitable splitting of the kernels for efficient discretisation to be applied. An advantage with the method to be presented is that the original problem is transformed into equations over the boundary of the domain $D$, thus reducing the dimensionality compared with discretizing the whole of $D$. When transforming with respect to time, one would expect the use of volume potentials for integral formulations of the stationary problem but this is circumvented as explained below. 

In \cite{Sun} the similar Cauchy problem studied, and a method is presented based on separation of variables in time and space, and further writing the spacewise part as a sum of two wave-like solutions. We do not use any such separation or splitting.

For the outline of the present work, in Section~2, we derive a boundary integral formulation of~(\ref{h_eq})--(\ref{h_bc}) in the form of a sequence of systems of boundary integral equations. This sequence is obtained by first applying the Laguerre transformation in time, and then for the obtained sequence of stationary problems invoke what is known as a fundamental sequence. Using this fundamental sequence, the solution to the stationary problems can be represented in terms of a sequence of single-layer potentials with densities to be determined, and no volume potential is needed which is an advantage. Matching the given transformed Cauchy data, a sequence of systems of boundary integral equations is obtained for the unknown densities. Analysis of the singularities of the involved kernels reveals that the mapping corresponding to this system is injective and have dense range, see~Theorem~\ref{op_prop}. As a preparation for the discretization, also included in Section~2 is a rewriting of the kernels that makes the singularities appear in explicit form. 

In Section~3, we introduce a parametrization of the boundary curves $\Gamma_1$ and $\Gamma_2$. With this parametrization, we can further rewrite the kernels and transform the system into $2\pi$-periodic boundary integral equations. 

Section~4 contains discretisation of the obtained integral equations, generating linear systems to solve for values of the densities. Tikhonov regularization is invoked for the stable solution of this system. Moreover, explicit formulas are given for the function values and traction on the inner boundary $\Gamma_1$.  

Numerical results are presented in Section~5, confirming that the outlined approach is a feasible way of reconstructing the missing data on $\Gamma_1$. We point out that some derivations of formulas in the present work are long and it is not possible, in order to keep the presentation at reasonable length, to give full details but we do refer to work where more details can be found. The exactness of the numerical results is a further confirmation of the correctness of the stated formulas. Some conclusions are stated in Section~6.

\section{Combination of the Laguerre transform and the boundary integral equation method}\label{Sect_2}
 
 We search for the solution $u$ of \eref{h_eq}--\eref{h_bc} as the (scaled) Fourier expansion with respect to the Laguerre polynomials, that is an expansion of the form 
\bql{lag6}
 u(x,t)=\kappa\sum_{n=0}^\infty
 u_n(x)L_n(\kappa t),
 \eq
 where 
 \bql{lag7}
  u_n(x)=\int_0^\infty e^{-\kappa t}L_n(\kappa t)u(x,t)\,dt,\quad 
n=0,1,2,\dots. 
 \eq
 Here, $L_n$ is the Laguerre polynomial of order $n$ with scaling parameter $\kappa>0$.

For the Fourier--Laguerre coefficients $u_n$ in~\eref{lag7} of the function $u$, using the recurrence relations for the Laguerre polynomials, it can be shown (see~\cite{ChMi}) that they satisfy the following sequence of Cauchy problems
\bql{s_e}
\Delta^* {u_n} - {\kappa^2}{u_n} = \sum_{m=0}^{n-1}\beta_{n-m}u_m , \quad \mbox{in } D,  
\eq
\bql{s_bc}
{u_n} = {f_{2,n}},\quad \text{on } \Gamma_2, \qquad Tu_n=g_{2,n}, \quad  \text{on } \Gamma_2, 
\eq
where $n = 0,1, \ldots$,   $\beta_{n}=\kappa^2 (n+1)$ and $\{f_{2,n}\}$ and $\{g_{2,n}\}$ are the Fourier-Laguerre sequences of coefficients of the given functions $f_2$ and $g_2$ in~\eref{h_bc}. It is straightforward to check the following result. 
\begin{theorem} 
A sufficiently smooth function~\eref{lag6} is the solution of the time-dependent problem \eref{h_eq}--\eref{h_bc} if and only if its Fourier-Laguerre coefficients  
$u_n$ for $n = 0,1,\ldots$, solve the sequence of stationary problems \eref{s_e}--\eref{s_bc}. 
\end{theorem}
Due to this, we focus on solving the sequence of stationary Cauchy problems~\eref{s_e}--\eref{s_bc}.

The work~\cite{ChaJoh} follows a string of works, where the solution to various Cauchy problems are represented in terms of a single-layer potential (see~\cite{ChaJoh16} for an overview). We use the same strategy here for~\eref{s_e}--\eref{s_bc}. As a suitable representation, we follow~\cite{ChMi} and shall therefore construct a solution to~\eref{s_e}--\eref{s_bc} in the form of a sequence of single-layer potentials
\bql{s_p}
{u_n}(x) = \frac{1}{{2\pi}}\sum_{\ell=1}^2\sum\limits_{m = 0}^n {\int\limits_{\Gamma_\ell} {{E _{n - m}}(x,y){q^\ell_m}(y)\,ds(y)}} ,
\quad
x \in D,
\eq
where $E_n$ is a sequence of fundamental solutions to~\eref{s_e}. We shall give an explicit expression for $E_n$ in Section~\ref{Sect_E_n}. In the next section, we derive a system for determining the densities $q^\ell_m$. Note that as a convention in this work, an upper index of a sub-indexed element does not refer to a power but is treated as an index unless otherwise explicitly stated.

\subsection{A system of boundary integral equations for~\eref{s_e}--\eref{s_bc}}
The boundary integral operators corresponding to the representation \eref{s_p} admit the same jump relations as the classical single-layer operator for the Laplace equation; this can be verified by noticing 
that each function in the fundamental sequence has at most a logarithmic singularity. We shall make a detailed analysis of the singularities in the following sections.
Matching~\eref{s_p} with the data~\eref{s_bc} and employing the jump properties, we obtain the following system of boundary integral equations 
\begin{equation}\label{sie}
\left\{\begin{array}{lr}
       \displaystyle{ \frac{1}{2\pi}\sum_{\ell=1}^2\int_{\Gamma_\ell}  E_0(x, y)q_n^\ell(y)\, ds(y)  = F_n (x)},  & x\in\Gamma_2,\\
     \displaystyle{   \frac{1}{2}q_n^2(x)+\frac{1}{2\pi}\sum_{\ell=1}^2\int_{\Gamma_\ell} T_xE_0(x, y)q_n^\ell(y) \, ds(y) = G_n (x)}, & x\in\Gamma_2 ,
        \end{array}\right.
\end{equation}
 for $n=0,\ldots,N$, with the right-hand sides
 $$
 \begin{array}{rcl}
F_n(x)& \ar =\ar & \displaystyle f_{2,n}(x)
 -\frac{1}{2\pi}\sum_{\ell=1}^2\sum_{m=0}^{n-1}\int_{\Gamma_\ell} E_{n-m}(x,y)q_m^\ell(y)\, ds(y)
 \end{array}
 $$
 and
 $$
  \begin{array}{rcl}
  G_n(x)& \ar =\ar & \displaystyle g_{2,n}(x)-\frac{1}{2}\sum_{m=0}^{n-1}q_m^2(x)
 -\frac{1}{2\pi}\sum_{\ell=1}^2\sum_{m=0}^{n-1}\int_{\Gamma_\ell} T_x E_{n-m}(x,y)q_m^\ell(y)\, ds(y).
 \end{array}
 $$
The unknown densities $q_m^1$ and $q_m^2$, $m=0,\ldots,N$, in \eref{s_p} are defined on the two (closed) boundary curves $\Gamma_1$ and $\Gamma_2$, respectively (note the above convention that upper indices do not denote a power).  The operator $T_x$ is as in~(\ref{stress}) with the sub-index indicating the variable with which derivatives shall be taken.

We have then reduced the Cauchy problem \eref{h_eq}--\eref{h_bc}, via the Laguerre transform~\eref{lag6} rendering the system~\eref{s_e}--\eref{s_bc} with solution in the form of the single-layer representation~\eref{s_p}, to the system of boundary integral equations~\eref{sie}. As mentioned above, in connection with jump properties, the kernels appearing in the integral equations in~\eref{sie} contain logarithmic singularities and, as we will see, $T_xE_n$ has in addition a strong singularity (Cauchy type); we take these singularities into account when proposing numerical discretization. 

Following the proof of~\cite[Theorem 4.1]{CaKr}, it is possible to verify that the corresponding operator matrix of the system \eref{sie}, built from the involved integral operators, has the following properties:
\begin{theorem}\label{op_prop}
The operator corresponding to the system\/~\eref{sie} is injective and has dense range, as a mapping between $L^2$-spaces on the boundary.
\end{theorem}
This result implies that Tikhonov regularization can be applied to solve \eref{sie} in a stable way.

\subsection{The sequence of fundamental solutions $E_n$ in~\eref{s_p}}\label{Sect_E_n}
The definition of $E_n$ is that it solves~(\ref{s_e}) but with a $2\times2$ diagonal matrix introduced in the right-hand side with the Dirac delta function $\delta(x-y)$ as diagonal elements of that matrix, see further~\cite[Definition~1]{ChMi}. The explicit expression for $E_n$ is given in~\cite[Theorem~3]{ChMi}, we recall that expression here. We have
\bql{E_n} 
E_n(x,y)=\Phi_{1,n}(|x-y|)I + \Phi_{2,n}(|x-y|)J(x-y),
\eq 
where $I$ is the identity matrix and
\bql{J}
J(x) = \frac{x \, x^\top}{|x|^2}, \quad x\in\real\setminus \{0\}.
\eq
Here,
\begin{equation}\label{Phi_l_n}
\begin{split}
\Phi_{\ell,n}(r) = {} & \frac{(-\ell)^{\ell-1}}{\kappa^2 r^2}\sum_{k=-2}^2 \chi_{k,n}\left(\Phi_{n+k}(\tfrac{\kappa}{c_s},r) - \Phi_{n+k}(\tfrac{\kappa}{c_p},r)\right)+\frac{(-1)^{\ell-1}}{c_p^2}\Phi_{n}(\tfrac{\kappa}{c_p},r) \\
& + \frac{\ell-1}{c_s^2}\Phi_{n}(\tfrac{\kappa}{c_s},r),
\end{split}
\end{equation}
for $\ell=1,2$, $\chi_{-2,n}=n(n-1)$, $\chi_{-1,n}=-4n^2$, $\chi_{0,n}=2(3n^2+3n+1)$, $\chi_{1,n}=-4(n+1)^2$ and $\chi_{2,n}=(n+1)(n+2)$
 and
\bql{Phi_n}
 \Phi_n(\gamma, r)=K_0(\gamma r)\,v_n(\gamma ,r)
 +K_1(\gamma r)\,w_n(\gamma, r),
\eq
where $K_0$ and $K_1$ are the modified Hankel functions of order zero and one, respectively (the modified Hankel function is sometimes termed differently such as Bessel function of the third kind, Basset's function or Macdonald's function), and 
\begin{equation}\label{pol}
 v_n(\gamma, r)=
 \sum_{m=0}^{\left[\frac{n}{2}\right]}a_{n,2m}(\gamma)r^{2m},\quad
 \quad
 w_n(\gamma, r)=
 \sum_{m=0}^{\left[\frac{n-1}{2}\right]}a_{n,2m+1}(\gamma)r^{2m+1},
\end{equation}
 \noindent
 for $n=0,1,\ldots,N-1$ ($w_0=0$)
 with the coefficients $a_{n,m}$ satisfying the recurrence relations
\begin{equation*}
\begin{aligned}
a_{n,0}(\gamma) &=1, & n &=0,1,\ldots,N-1, \\
a_{n,n}(\gamma) &=-\frac{\gamma}{n}\;a_{n-1,n-1}(\gamma),&  n&=1,2,\ldots,N-1
\end{aligned}
\end{equation*}
and 
\bql{a_nm}
 a_{n,m}(\gamma)=\frac{1}{2\gamma m}
 \left\{4\left[\frac{m+1}{2}\right]^2a_{n,m+1}(\gamma)
 -\gamma^2\sum_{k=m-1}^{n-1}(n-k+1) a_{k,m-1}(\gamma)\right\}.
\eq

The above expressions are derived in~\cite{ChMi} by calculating the Laguerre transform of the fundamental solution to the elastodynamic equation~(\ref{h_eq}). A key part is that~(\ref{Phi_n}) is a fundamental sequence for the Laguerre transformation of the wave equation; derivations that lead up to expressions of the form~(\ref{Phi_n}) for fundamental sequences of Laguerre time-transformed equations are given in~\cite{ChKr}.

\subsection{Calculation of the traction of the fundamental sequence $E_n$}
In \eref{sie}, the traction of the fundamental sequence $E_n$ is needed. We derive an explicit expression for it. We note first that for a function $f$ and a matrix $G$  there is the product rule (see \cite[Sect.~2]{ChKrMo})
\bql{stress1}
T(fG)=T(fI)\,G+f\,TG,
\eq
where, for a matrix, the convention is that the traction operator $T$ from~\eref{stress} acts on each column of that matrix. We recall from~\cite[Sect.~2]{ChKrMo} when $f:(0,\infty)\to\real$ is continuously differentiable, a
lengthy but straightforward calculation (employing~\eref{stress1}) reveals
\bql{stress2}
T_x (f(|x-y|)\,I )=
\frac{f^\prime(|x-y|)}{|x-y|}\;U_1(x,y),
\eq
where 
\bql{stress3}
U_1(x,y)=
\lambda \nu(x)\,(x-y)^\top
 + \mu (x-y)\,\nu(x)^\top +\mu \nu(x)^\top\,(x-y)\,I.
\eq
Following~\cite{ChKrMo}, we also have 
\bql{stress4}
T_x (J(x-y))=
\frac{1}{|x-y|^2}\;
U_2(x,y),
\eq
where $J$ is as in~\eref{J} and
\bql{stress5}
U_2(x,y)=
(\lambda +2\mu)\,\nu(x)\,(x-y)^\top
+\mu\, (x-y)\,\nu(x)^\top+\mu\,\nu(x)^\top\,(x-y)[I-4J(x-y)].
\eq

Using the expression~\eref{E_n} for $E_n$ and the definition of the traction operator~\eref{stress}, we find by applying~\eref{stress2} and~\eref{stress4} that the sought after expression for the traction of $E_n$ is
\bql{TE}
T_xE_{n}(x,y)=\frac{U_1(x,y)}{|x-y|^2}\left[\tilde{\Phi}_{1,n}(|x-y|)I+\tilde{\Phi}_{2,n}(|x-y|)J(x-y)\right]+
\Phi_{2,n}(|x-y|)\frac{U_2(x,y)}{|x-y|^2}
\eq
with
\begin{equation}\label{Phit_e_n}
\begin{split}
\tilde{\Phi}_{\ell,n}(r) = {} &\frac{(-\ell)^{\ell-1}}{\kappa^2 r^2}\sum_{k=-2}^2 \chi_{k,n}\left(\tilde{\Phi}_{n+k}(\tfrac{\kappa}{c_s},r) -2\Phi_{n+k}(\tfrac{\kappa}{c_s},r)+ 2\Phi_{n+k}(\tfrac{\kappa}{c_p},r) -\tilde{\Phi}_{n+k}(\tfrac{\kappa}{c_p},r)\right)\\
& +\frac{(-1)^{\ell-1}}{c_p^2}\tilde{\Phi}_{n}(\tfrac{\kappa}{c_p},r) +\frac{\ell-1}{c_s^2}\tilde{\Phi}_{n}(\tfrac{\kappa}{c_s},r),
\end{split}
\end{equation}
for $\ell=1,2$ and
\begin{equation}\label{Phi_n_t}
 \tilde{\Phi}_n(\gamma, r)=K_0(\gamma r)\,\tilde{v}_n(\gamma ,r)
 +K_1(\gamma r)\,\tilde{w}_n(\gamma, r).
\end{equation}
Here, we introduced polynomials
\begin{equation}\label{polt}
\begin{split}
  \tilde{v}_n(\gamma,r) &=2 \sum_{m=1}^{\left[\frac{n}{2}\right]}ma_{n,2m}(\gamma)r^{2m} -\gamma \sum_{m=0}^{\left[\frac{n-1}{2}\right]}a_{n,2m+1}(\gamma)r^{2m+2}, \\
   \tilde{w}_n(\gamma,r) &=2 \sum_{m=1}^{\left[\frac{n-1}{2}\right]}ma_{n,2m+1}(\gamma)r^{2m+1}-\gamma \sum_{m=0}^{\left[\frac{n}{2}\right]}a_{n,2m}(\gamma)r^{2m+1}.
\end{split}
\end{equation}
The expression for $\tilde{\Phi}_{\ell,n}$, with accompanying function~(\ref{Phi_n_t}) and~polynomials~(\ref{polt}), is obtained by simply differentiating~\eref{Phi_l_n} using the expressions~(\ref{Phi_n})--(\ref{pol}) and the relation for derivatives of the modified Hankel functions~\cite[Eqns. 9.6.27--28]{AbSt}.

\subsection{Analysis of singularities of the kernels in the system~\eref{sie}}\label{Sect_2_4}
The modified Hankel functions have the following series representations (\cite[Eqns. 9.6.11, 9.6.13]{AbSt})
 \begin{equation}\label{Macdonald}
  K_0(z) = - \left(\ln \frac{z}2+C\right)\,I_0(z)
 + S_0(z),
\quad
 K_1(z)= \frac{1}{z}+\left(\ln \frac{z}2 +C\right)\,I_1(z)
+S_1(z),
 \end{equation}
with
$$
I_0(z)=\sum^\infty_{n=0} \;
 \frac{1}{(n!)^2}\,\left(\frac{z}{2}\right)^{2n},
 \quad
 I_1(z)=\sum^\infty_{n=0} \;
 \frac{1}{n!(n+1)!}\,\left(\frac{z}{2}\right)^{2n+1},
 $$
and
$$
S_0(z)=\sum^\infty_{n=1}
 \frac{\psi(n)}{ (n!)^2} \,\left(\frac{z}{2}\right)^{2n}, \quad
S_1(z)=-\frac{1}{2}\sum^\infty_{n=0}
 \frac{\psi(n+1)+\psi(n)}{ n!(n+1)!}\,\left(\frac{z}{2}\right)^{2n+1}.
$$
Here, we put $\psi(0)=0$,
 $$
 \psi(n)=\sum_{m=1}^n\frac{1}{m}\;, \quad n=1,2,\ldots,
 $$
 and let $C = 0.57721\ldots$ denote Euler's constant. Thus, using~(\ref{Macdonald}) we can rewrite the functions $\Phi_n$ in~\eref{Phi_n} as
\begin{equation}\label{Phi_n_r}
\Phi_n(\gamma, r)= \phi_n(\gamma,r)\ln r+\varphi_n(\gamma,r), \quad n=0,1,\ldots,
\end{equation}
where
\begin{equation}\label{phi}
\phi_n(\gamma,r)=-I_0(\gamma r)v_n(\gamma, r)+I_1(\gamma r)w_n(\gamma, r)
\end{equation}
and 
\begin{equation}\label{varphi_n}
\begin{split}
\varphi_n(\gamma, r) = { } & \Big[-\Big(C+\ln\frac{\gamma}{2}\Big)I_0(\gamma r)+S_0(\gamma r)\Big]v_n(\gamma, r)\\
& +\Big[\frac{1}{\gamma r}+\Big(C+\ln\frac{\gamma}{2}\Big)I_1(\gamma r)+S_1(\gamma r)\Big]w_n(\gamma, r).
\end{split}
\end{equation}

The representation~(\ref{Phi_n_r}) of $\Phi_n$ then implies that we have the following expressions for the functions defined in~\eref{Phi_l_n},
\bql{Phi_ln}
\Phi_{\ell,n}(r)=\eta_{\ell,n}(r)\ln r +\xi_{\ell,n}(r),\quad \ell=1,2
\eq
with
\begin{equation}\label{eta_l_n}
\begin{split}
\eta_{\ell,n}(r) = {} & \frac{(-\ell)^{\ell-1}}{\kappa^2 r^2}\sum_{k=-2}^2 \chi_{k,n}\left(\phi_{n+k}(\tfrac{\kappa}{c_s},r) - \phi_{n+k}(\tfrac{\kappa}{c_p},r)\right)+
\frac{(-1)^{\ell-1}}{c_p^2}\phi_{n}(\tfrac{\kappa}{c_p},r) \\
& +\frac{\ell-1}{c_s^2}\phi_{n}(\tfrac{\kappa}{c_s},r) ,
\end{split}
\end{equation}
and
\begin{equation}\label{xi_l_n}
\begin{split}
\xi_{\ell,n}(r) = {} & \frac{(-\ell)^{\ell-1}}{\kappa^2 r^2}\sum_{k=-2}^2 \chi_{k,n}\left(\varphi_{n+k}(\tfrac{\kappa}{c_s},r) - \varphi_{n+k}(\tfrac{\kappa}{c_p},r)\right)+
\frac{(-1)^{\ell-1}}{c_p^2}\varphi_{n}(\tfrac{\kappa}{c_p},r)\\
& 
+\frac{\ell-1}{c_s^2}\varphi_{n}(\tfrac{\kappa}{c_s},r).
\end{split}
\end{equation}
Analysis of the representations \eref{eta_l_n} and \eref{xi_l_n} shows that we have the following asymptotic behavior with respect to $r$ (see \cite{ChMi})
$$
\eta_{\ell,n}(r)=\eta_{\ell,n}(0)+\mathcal{O}(r^2), \qquad \xi_{\ell,n}(r)=\xi_{\ell,n}(0) +\mathcal{O}(r^2)
$$
with
\begin{align*}
\eta_{\ell,n}(0)&=\frac{(-\ell)^{\ell-1}}{\kappa^2}\sum_{k=-2}^2 \chi_{k,n}\left(\epsilon_{n+k,2}(\tfrac{\kappa}{c_s})-\epsilon_{n+k,2}(\tfrac{\kappa}{c_p})\right)+\frac{(-1)^{\ell-1}}{c_p^2}\epsilon_{n,0}(\tfrac{\kappa}{c_p})
\\
&\phantom{=}+\frac{\ell-1}{c_s^2}\epsilon_{n,0}(\tfrac{\kappa}{c_s})
\end{align*}
and
\begin{align*}
\xi_{\ell,n}(0)&=\frac{(-\ell)^{\ell-1}}{\kappa^2}\sum_{k=-2}^2 \chi_{k,n}\left(\varepsilon_{n+k,2}(\tfrac{\kappa}{c_s})-\varepsilon_{n+k,2}(\tfrac{\kappa}{c_p})\right)+\frac{(-1)^{\ell-1}}{c_p^2}\varepsilon_{n,0}(\tfrac{\kappa}{c_p})
\\
&\phantom{=}+\frac{\ell-1}{c_s^2}\varepsilon_{n,0}(\tfrac{\kappa}{c_s}).
\end{align*}
Here
$$
\epsilon_{n,0}(\gamma)=-a_{n,0}(\gamma), \quad 
\epsilon_{n,2}(\gamma)=-\frac{\gamma^2}{4}a_{n,0}(\gamma)+\frac{\gamma}{2}a_{n,1}(\gamma)-a_{n,2}(\gamma)
$$
and
\begin{align*}
\varepsilon_{n,0}(\gamma) &=-\Big(C+\ln\frac{\gamma}{2}\Big)a_{n,0}(\gamma)+\frac{1}{\gamma}a_{n,1}(\gamma), \\
\varepsilon_{n,2}(\gamma) &=\Big(C+\ln\frac{\gamma}{2}\Big)\Big(-\frac{\gamma^2}{4}a_{n,0}(\gamma)+\frac{\gamma}{2}a_{n,1}(\gamma)-a_{n,2}(\gamma)\Big)
+\frac{\gamma^2}{4}a_{n,0}\\
&\phantom{=}-\frac{\gamma}{4}a_{n,1}(\gamma)+\frac{1}{\gamma}a_{n,3}(\gamma).
\end{align*}

Note here that a straightforward calculation using the recurrence  formula~\eref{a_nm} gives $\eta_{1,n}(0)=-1/(2 c_p^2) - 1/(2 c_s^2)$, $\eta_{2,n}(0)=0$ and $\xi_{2,n}(0)=-1/(2 c_p^2) + 1/(2 c_s^2)$ (here and below the upper index for $c_p$ and $c_s$ denote of course a power).

Using~\eref{Phi_ln} in the definition of the fundamental sequence~\eref{E_n}, this fundamental sequence can be written as 
\begin{equation}\label{E_n_form}
\begin{split}
E_{n}(x,y)= {} & \ln |x-y|\Big[\eta_{1,n}(|x-y|)I+\eta_{2,n}(|x-y|)J(x-y)\Big]\\
&
+\xi_{1,n}(|x-y|)I+\xi_{2,n}(|x-y|)J(x-y).
\end{split}
\end{equation}
Thus, we have verified that the sequence $E_n$ has a singularity of logarithmic type.

Turning to the traction, the similar analysis applied to~\eref{TE}--\eref{polt}, reveals that~(\ref{Phi_n_t}) can be written 
\begin{equation}\label{Phi_n_t_r}
\tilde{\Phi}_n(\gamma, r) = \tilde{\phi}_n(\gamma,r)\ln r+\tilde{\varphi}_n(\gamma,r), \quad n=0,1,\ldots,
\end{equation}
where
\begin{equation}\label{phi_t}
\tilde{\phi}_n(\gamma,r)=-I_0(\gamma r)\tilde{v}_n(\gamma, r)+I_1(\gamma r)\tilde{w}_n(\gamma, r)
\end{equation}
and 
\begin{equation}\label{varphi_n_t}
\begin{split}
\tilde{\varphi}_n(\gamma, r)= {} & \Big[-\Big(C+\ln\frac{\gamma}{2}\Big)I_0(\gamma r)+S_0(\gamma r)\Big]\tilde{v}_n(\gamma, r)\\
& +\Big[\frac{1}{\gamma r}+\Big(C+\ln\frac{\gamma}{2}\Big)I_1(\gamma r)+S_1(\gamma r)\Big]\tilde{w}_n(\gamma, r).
\end{split}
\end{equation}
Using~\eref{Phi_n_t_r} in~\eref{Phit_e_n}, we can then derive the following representation for the functions $\tilde{\Phi}_{\ell,n}$,
$$
\tilde{\Phi}_{\ell,n}(r)=\tilde{\eta}_{\ell,n}(r)\ln r +\tilde{\xi}_{\ell,n}(r),\quad \ell=1,2
$$
with
\begin{align*}
\tilde{\eta}_{\ell,n}(r) &=\frac{(-\ell)^{\ell-1}}{\kappa^2 r^2}\sum_{k=-2}^2 \chi_{k,n}\left(\tilde{\phi}_{n+k}(\tfrac{\kappa}{c_s},r) 
-2\phi_{n+k}(\tfrac{\kappa}{c_s},r) + 2\phi_{n+k}(\tfrac{\kappa}{c_p},r)-\tilde{\phi}_{n+k}(\tfrac{\kappa}{c_p},r)\right) \\
&\phantom{=}+
\frac{(-1)^{\ell-1}}{c_p^2}\tilde{\phi}_{n}(\tfrac{\kappa}{c_p},r)+\frac{\ell-1}{c_s^2}\tilde{\phi}_{n}(\tfrac{\kappa}{c_s},r) , \\
\tilde{\xi}_{\ell,n}(r) &=\frac{(-\ell)^{\ell-1}}{\kappa^2 r^2}\sum_{k=-2}^2 \chi_{k,n}\left(\tilde{\varphi}_{n+k}(\tfrac{\kappa}{c_s},r) 
-2\varphi_{n+k}(\tfrac{\kappa}{c_s},r) + 2\varphi_{n+k}(\tfrac{\kappa}{c_p},r)-\tilde{\varphi}_{n+k}(\tfrac{\kappa}{c_p},r)\right) \\
&\phantom{=}+
\frac{(-1)^{\ell-1}}{c_p^2}\tilde{\varphi}_{n}(\tfrac{\kappa}{c_p},r)+\frac{\ell-1}{c_s^2}\tilde{\varphi}_{n}(\tfrac{\kappa}{c_s},r),
\end{align*}
where $\phi_n$ is defined in~(\ref{phi}), $\tilde{\phi}_n$ in~(\ref{phi_t}), $\varphi_n$ in~(\ref{varphi_n}) and $\tilde{\varphi}_n$ in~(\ref{varphi_n_t}).

It is straightforward to see that 
$$
\tilde{\xi}_{\ell,n}(r)=\tilde{\xi}_{\ell,n}(0)+\mathcal{O}(r^2), \quad \tilde{\eta}_{\ell,n}(r)=\mathcal{O}(r^2),
$$
where
\begin{align*}
\tilde{\xi}_{\ell,n}(0)&=\frac{(-\ell)^{\ell-1}}{\kappa^2}\sum_{k=-2}^2 \chi_{k,n}\left(\tilde{\varepsilon}_{n+k,2}(\tfrac{\kappa}{c_s})-\tilde{\varepsilon}_{n+k,2}(\tfrac{\kappa}{c_p})+2[\varepsilon_{n+k,2}(\tfrac{\kappa}{c_s})-\varepsilon_{n+k,2}(\tfrac{\kappa}{c_p})]\right)
\\
&\phantom{=}+\frac{(-1)^{\ell-1}}{c_p^2}\tilde\varepsilon_{n,0}(\tfrac{\kappa}{c_p}) +\frac{\ell-1}{c_s^2}\tilde\varepsilon_{n,0}(\tfrac{\kappa}{c_s}),
\end{align*}
with
\begin{align*}
\tilde\varepsilon_{n,0}(\gamma) &= -a_{n,0}, \\
\tilde{\varepsilon}_{n,2}(\gamma) &=\left(C+\ln\frac{\gamma}{2}\right)\left[\gamma a_{n,1}(\gamma)-2a_{n,2}(\gamma)-\frac{\gamma^2}{2}a_{n,0}(\gamma)\right] + \frac{\gamma^2}{4}a_{n,0}(\gamma) \\
&\phantom{=}- a_{n,2}(\gamma) +\frac{2}{\gamma} a_{n,3}(\gamma).
\end{align*}

Again, a straightforward calculation shows that $\tilde{\xi}_{2,n}(0)=0$ and $\tilde{\xi}_{1,n}(0)=-1/(2 c_p^2) - 1/(2 c_s^2)$.
Thus, we have arrived at the final form of the traction of the fundamental sequence, which we shall use in subsequent sections, and it can be stated as 
\begin{equation}\label{TE_form}
T_xE_{n}(x,y)=  \ln|x-y|W^{n,1}(x,y)+W^{n,2}(x,y)
\end{equation}
with matrices
\begin{align*}
W^{n,1}(x,y)&=\frac{U_1(x,y)}{|x-y|^2}\left[\tilde{\eta}_{1,n}(|x-y|)I+\tilde{\eta}_{2,n}(|x-y|)J(x-y)\right]+\eta_{2,n}(|x-y|)\frac{U_2(x,y)}{|x-y|^2}, \\
W^{n,2}(x,y)&=\frac{U_1(x,y)}{|x-y|^2}\left[\tilde{\xi}_{1,n}(|x-y|)I+\tilde{\xi}_{2,n}(|x-y|)J(x-y)\right]+\xi_{2,n}(|x-y|)\frac{U_2(x,y)}{|x-y|^2}
\end{align*}
and~$U_1$ and $U_2$ defined in~(\ref{stress3}) and~(\ref{stress5}), respectively. Thus, we have verified that the elements in the sequence $T_xE_n$ has a logarithmic singularity and a strong singularity. 

\section{Transformation to 2$\pi$-periodic integral equations}
In order to apply standard quadrature rules~\cite{Kr} for singular periodic integrals, we introduce a suitable parametrization of the boundary. The system~\eref{sie} is then rewritten using this parameterisation, and we state expressions for the parametrized kernels.

We assume that the boundary curves $\Gamma_\ell$, $\ell=1,2$, are  sufficiently smooth and given by a parametric representation
$$
\Gamma_\ell=\{x_\ell(s)=(x_{1\ell}(s),x_{2\ell}(s)):\, s\in[0,2\pi]\}.
$$  
The system~\eref{sie} can then be written in parametric form
\begin{equation}\label{psie}
\left\{\begin{array}{lr}
       \displaystyle{ \frac{1}{2\pi}\sum_{\ell=1}^2\int_0^{2\pi}  H_{\ell,2}^0(s,\sigma)\psi_n^\ell(\sigma)\, d\sigma  = F_n (s)},  & s\in [0,2\pi],\\
     \displaystyle{   \frac{\psi_n^2(s)}{2|x'_2(s)|}+\frac{1}{2\pi}\sum_{\ell=1}^2\int_0^{2\pi}  Q_{\ell,2}^0(s,\sigma)\psi_n^\ell(\sigma)\, d\sigma = G_n (s)}, & s\in [0,2\pi] ,
        \end{array}\right.
\end{equation}
 for $n=0,\ldots,N$, where $\psi_n^\ell(s)=q_n^\ell(x_\ell(s))|x'_\ell(s)|$.
 
The right-hand sides in~(\ref{psie}) are given by
 $$
 \begin{array}{rcl}
F_n(s)& \ar =\ar & \displaystyle f_{2,n}(x_2(s))
 -\frac{1}{2\pi}\sum_{\ell=1}^2\sum_{m=0}^{n-1}\int_0^{2\pi} H_{\ell,2}^{n-m}(s,\sigma)\psi_m^\ell(\sigma)\, d\sigma
 \end{array}
 $$
 and
 $$
  \begin{array}{rcl}
  G_n(s)& \ar =\ar & \displaystyle g_{2,n}(x_2(s))-\frac{1}{2|x_2'(s)|}\sum_{m=0}^{n-1}\psi_m^2(s)
-\frac{1}{2\pi}\sum_{\ell=1}^2\sum_{m=0}^{n-1}\int_0^{2\pi} Q_{\ell,2}^{n-m}(s,\sigma)\psi_m^\ell(\sigma)\, d\sigma.
 \end{array}
 $$
 The kernels in~\eref{psie} are
 \begin{equation}\label{def_par_ker}
 H_{\ell,k}^n(s,\sigma)=E_n (x_k(s),x_\ell(\sigma)) 
 \end{equation} 
 and
 \begin{equation}\label{def_par_ker_Q} 
 Q_{\ell,k}^n(s,\sigma)= T_x E_n(x_k(s),x_\ell(\sigma)) ,
 \end{equation}
 for $s\neq \sigma$, $\ell,k=1,2$, $n=0,\ldots,N$, and $E_n$ the fundamental sequence~(\ref{E_n}) and the traction $T_x$ given by~(\ref{stress}).

\subsection{Expressions for the singular kernel~$H^n_{\ell,\ell}$}
According to the analysis of the singularities undertaken in Section~\ref{Sect_2}, the fundamental sequence~(\ref{E_n}) can equivalently be written as in~(\ref{E_n_form}). Using this latter expression in combination with~(\ref{def_par_ker}), we can write
\begin{equation}\label{Hnellell}
   H^n_{\ell,\ell}(s, \sigma ) =
 \ln\left(  \frac{4}e \sin^2 \frac{s- \sigma}{2}\right)
 H^{n,1}_{\ell,\ell} (s, \sigma )
  + H^{n,2}_{\ell,\ell} (s,\sigma ),
\end{equation}
where
$$
H^{n,1}_{\ell,\ell} (s, \sigma )=\frac{1}{2}\left[\eta_{1,n}(|x_\ell (s)-x_\ell(\sigma)|)I+\eta_{2,n}(|x_\ell(s)-x_\ell(\sigma)|)J(x_\ell(s)-x_\ell(\sigma))\right]
$$
and
$$
H^{n,2}_{\ell,\ell} (s, \sigma )= H^{n}_{\ell,\ell} (s, \sigma )-
 \ln\left(  \frac{4}e \sin^2 \frac{s- \sigma}{2}\right)
 H^{n,1}_{\ell,\ell} (s, \sigma ),
$$
with the diagonal terms
$$
H^{n,1}_{\ell,\ell} (s, s)=\frac{1}{2}\eta_{1,n}(0)I
$$
and
$$
H^{n,2}_{\ell,\ell} (s, s )= \frac12 \ln(|x_\ell'(s)|^2 e)\eta_{1,n}(0)I+\xi_{1,n}(0)I+\xi_{2,n}(0)\widetilde J_\ell (s).
$$
Here, we used the Taylor expansion of the matrix $J$ from~\eref{J} resulting in the term
$$
\widetilde J_\ell (s)=
\frac{x_\ell^\prime(s) x_\ell^\prime(s)^\top}{|x_\ell^\prime(s)|^2}.
$$

\subsection{Expressions for the singular kernel~$Q_{\ell,\ell}^n$}
In Section~\ref{Sect_2}, it was shown that the expression~\eref{TE} for the traction of the fundamental sequence can be equivalently written as~\eref{TE_form}. From this latter expression it follows that the kernels $Q_{\ell,\ell}^n$ defined in~\eref{def_par_ker_Q} have logarithmic and strong type singularities.  To handle those we note the following expansions
$$
\frac{U_k(x_\ell(s),x_\ell(\sigma))}{|x_\ell(s)-x_\ell(\sigma)|^2}=\frac{1}{|x_\ell'(s)|^2(s-\sigma)}\tilde{U}_{k,\ell}(s)-
\frac{1}{2|x_\ell'(s)|^2}\hat{U}_{k,\ell}(s)+\frac{x_\ell'(s)\cdot x_\ell''(s)}{|x_\ell'(s)|^4}\tilde{U}_{k,\ell}(s)+\mathcal{O}(s-\sigma),
$$
where
\begin{align*}
\tilde{U}_{1,\ell}(s) &=\lambda \nu(x_\ell(s))x'_\ell(s)^\top+\mu x'_\ell(s)\nu(x_\ell(s))^\top, \\
\tilde{U}_{2,\ell}(s) &=(\lambda+2\mu) \nu(x_\ell(s))x'_\ell(s)^\top+\mu x'_\ell(s)\nu(x_\ell(s))^\top,
\end{align*}
and
\begin{align*}
\hat{U}_{1,\ell}(s) &=\lambda \nu(x_\ell(s))x''_\ell(s)^\top+\mu x''_\ell(s)\nu(x_\ell(s))^\top+\mu\nu(x_\ell(s))^\top x''_\ell(s)I , \\
\hat{U}_{2,\ell}(s) &=(\lambda+2\mu) \nu(x_\ell(s))x''_\ell(s)^\top+\mu x''_\ell(s)\nu(x_\ell(s))^\top+\mu\nu(x_\ell(s))^\top x''_\ell(s)(I-4\tilde{J}_\ell (s)).
\end{align*}

We therefore write
%
%
\begin{equation}\label{Qnellell}
 Q_{\ell,\ell}^n(s,\sigma)= \ln \left( \frac{4}{e}\sin^2 \frac{s-\sigma}{2}\right) Q_{\ell,\ell}^{n,1}(s,\sigma)+\cot\frac{\sigma-s}{2} Q_{\ell,\ell}^{n,2}(s)+Q_{\ell,\ell}^{n,3}(s,\sigma),
\end{equation}
where
\begin{align*}
Q_{\ell,\ell}^{n,1}(s,\sigma) &=\frac{1}{2}W^{n,1}(x_\ell(s),x_\ell(\sigma)), \\
Q_{\ell,\ell}^{n,2}(s) &=-\frac{1}{2|x'_\ell(s)|^2}\left[\tilde{\xi}_{1,n}(0)\tilde{U}_{1,\ell}(s)+\xi_{2,n}(0)\tilde{U}_{2,\ell}(s)\right],
\end{align*}
and
$$
 Q_{\ell,\ell}^{n,3}(s,\sigma)=Q_{\ell,\ell}^n(s,\sigma)-\ln \left( \frac{4}{e}\sin^2 \frac{s-\sigma}{2}\right) Q_{\ell,\ell}^{n,1}(s,\sigma)-\cot\frac{\sigma-s}{2} Q_{\ell,\ell}^{n,2}(s).
$$ 
The  diagonal terms have the form
$$
Q_{\ell,\ell}^{n,1}(s,s)=0,
$$
and
\begin{align*}
Q_{\ell,\ell}^{n,3}(s,s)&=-\frac{1}{2|x'_\ell(s)|^2}\left(\tilde{\xi}_{1,n}(0)[\hat{U}_{1,\ell}(s)-\frac{2x'_\ell(s)\cdot x''_\ell(s)}{|x'_\ell(s)|^2}\tilde{U}_{1,\ell}(s)]\right.\\
&\left.\phantom{=}+\xi_{2,n}(0)[\hat{U}_{2,\ell}(s)-\frac{2x'_\ell(s)\cdot x''_\ell(s)}{|x'_\ell(s)|^2}\tilde{U}_{2,\ell}(s)]\right).
\end{align*}


\section{Full discretization} 
The effort in rewriting the kernels in~(\ref{psie}) using~(\ref{Hnellell}) and~(\ref{Qnellell}) now pays off in that we can employ the following standard quadrature rules~\cite{Kr} for numerical discretisation, 
\begin{align*}
 \frac{1}{2\pi}
 \int_0^{2\pi}
 f(\sigma)\,
 d\sigma
& \approx 
\frac1{2M} \sum_{k=0}^{2M-1}
  f(s_k), \\ 
 \frac{1}{2\pi}
 \int_0^{2\pi}
 f(\sigma)
 \ln\left(\frac{4}{e} \sin^2 {s- \sigma \over 2}\right)
 d\sigma
& \approx
 \sum_{k=0}^{2M-1}
 {R}_{k}(s)\,f(s_k)
\end{align*}
and
$$
 \frac{1}{2\pi}
 \int_0^{2\pi}
 f(\sigma)
 \cot\frac{\sigma - s}{2}\,
 d\sigma
 \approx
 \sum_{k=0}^{2M-1}
 {S}_{k}(s)\,f(s_k)
 $$
with mesh points 
\begin{equation}\label{mesh_points}
s_k=kh, \quad k=0,\ldots,2M-1, \quad h=\pi/M,
\end{equation} 
and the weight functions 
$$
 {R}_{k}(s) = 
\displaystyle
 - {1 \over 2M} \;  \left(1+2\sum^{M-1}_{m=1} \; {1 \over m} \,
\cos m(s- s_k) - {1 \over M} \, \cos M(s-s_k )\right),
$$
and
$$
S_k(s)={1 \over 2M}[1-(-1)^k \cos Ms)]\cot\frac{s_k - s}{2}, \quad s\ne s_k
$$
in order to approximate the boundary integrals in \eref{psie}. 

Collocating the approximation at the nodal points using the mesh points $\{s_k\}$ from~(\ref{mesh_points}) leads to the sequence of linear systems
 \bql{sle}
 \left\{
\begin{aligned}
&\sum\limits_{j = 0}^{2M-1}\left\{\frac{1}{2M} H_{1,2}^0(s_i, s_j)\psi_{n,j}^1 + 
	\left[R_j(s_i){H}_{2,2}^{0,1}(s_i, s_j)+\frac{1}{2M}{H}_{2,2}^{0,2}(s_i, s_j)\right]\psi_{n,j}^2\right\} = \tilde{F}_{n,i}, \\
&\sum\limits_{j = 0}^{2M-1} \left\{\frac{1}{2M}Q_{1,2}^0(s_i, s_j)\psi_{n,j}^1  + 
	\left[R_j(s_i){Q}_{2,2}^{0,1}(s_i, s_j)+S_j(s_i){Q}_{2,2}^{0,2}(s_i)\right.\right.\\
	&\left.\left.+\frac{1}{2M}{Q}_{2,2}^{0,3}(s_i, s_j)\right]\psi_{n,j}^2\right\}
	+ \frac{\psi_{n,i}^2}{2|x_2'(s_i)|}= \tilde{G}_{n,i}
\end{aligned}
\right.
 \eq
  for the unknown coefficients $\psi_{n,i}^l \approx \psi_n^l (s_i),\, \, i=0,\ldots,2M-1$, with the right-hand sides
\begin{equation}\label{Fni}
\begin{split}
\tilde{F}_{n,i}  = {} & f_{2,n}(x_2(s_i))
 -\sum_{j=0}^{2M-1}\sum_{m=0}^{n-1}\left\{ \frac{1}{2M}H_{1,2}^{n-m}(s_i, s_j)\psi_{m,j}^1 \right.\\
& \left. + \left[R_j(s_i){H}_{2,2}^{n-m,1}(s_i, s_j)+\frac{1}{2M}{H}_{2,2}^{n-m,2}(s_i, s_j)\right]\psi_{m,j}^2\right\}
\end{split} 
\end{equation}
 and
\begin{equation}\label{Gni}
\begin{split}
\tilde{G}_{n,i} = {} & g_{2,n}(x_2(s_i))-\frac{1}{2|x_2'(s_i)|}\sum_{m=0}^{n-1}\psi_{m,i}^2
 -\sum_{j=0}^{2M-1}\sum_{m=0}^{n-1}\left\{ \frac{1}{2M}Q_{1,2}^{n-m}(s_i, s_j)\psi_{m,j}^1\right.\\ 
& \left.+  \left[R_j(s_i){Q}_{2,2}^{n-m,1}(s_i, s_j)+S_j(s_i){Q}_{2,2}^{n-m,2}(s_i)+\frac{1}{2M}{Q}_{2,2}^{n-m,3}(s_i, s_j)\right]\psi_{m,j}^2\right\},
\end{split} 
\end{equation}
for $n=0,\ldots,N$.

 
We point out that the matrix corresponding to \eref{sle} is the same for every $n$ but having a recurrence right-hand side given by (\ref{Fni}) and (\ref{Gni}) containing the solutions of the previous systems. Clearly, each of the linear systems has a high-condition number since the Cauchy problem~\eref{s_e}--\eref{s_bc} is ill-posed and therefore Tikhonov regularization has to be incorporated.
 
We end this section by giving formulas for the Cauchy data on the boundary $\Gamma_1$ of the solution $u$ to~\eref{s_e}--\eref{s_bc}. Using~\eref{s_p} we have the following representation of the function value
 \begin{equation}\label{f1}
 f_{1,n}(x)=u_n(x)=\frac{1}{2\pi}\sum_{\ell=1}^2\sum_{m=0}^n
 \int_{\Gamma_\ell} E_{n-m}(x,y)q_m^\ell(y) \, ds(y),
 \quad x\in \Gamma_1
\end{equation}
 and for the traction
\begin{equation}\label{g1}
 g_{1,n}(x)=Tu_n(x)= -\frac{1}{2}\sum_{m=0}^{n}q_m^1(x)
 +\frac{1}{2\pi}\sum_{\ell=1}^2\sum_{m=0}^{n}\int_{\Gamma_\ell}T_xE_{n-m}(x,y)q_m^\ell(y)  ds(y), \quad x\in \Gamma_1.
\end{equation}
The numerical approximation of these expressions  
can be obtained using the given quadrature rules via similar calculations as those given above, with values of the densities generated from~(\ref{sle}).

 \section{Numerical examples}
\noindent
We present some numerical examples with the proposed method. The elastodynamic problem under consideration has a vector-valued function as solution. We shall not overload this section with results and figures for every component and domain considered but present enough data to have an idea about the accuracy and to have results to test against in case any reader wish to  implement the method. Further to this, the results given are not the best possible out of all tests done but rather typical. Thus, if the reader implements the method for corresponding configurations similar accuracy is expected. 
 
A key formula in our derivations, which is rather lengthy to derive, is~\eref{TE_form}. We shall therefore first start with a demonstration that  numerically verifies this formula. Thus, we consider the following sequence of well-posed stationary problems 
\begin{subequations}
\begin{alignat}{2}
\Delta^* {u_n} - {\kappa^2}{u_n} &= \sum_{m=0}^{n-1}\beta_{n-m}u_m , \quad &&\mbox{in }  D,  \label{s_e2}\\
Tu_n &= g_{1,n}, \quad &&
 \text{on } \Gamma_1, \label{s_bc2}\\
Tu_n &= g_{2,n}, \quad &&
 \text{on } \Gamma_2, \label{s_bc3}
\end{alignat}
\end{subequations}
where $g_{1,n} (x)= [TE_n (x,z_1)]_1, \, x\in \Gamma_1$ and $g_{2,n} (x)= [TE_n (x,z_1)]_1, \, x\in \Gamma_2$, for an arbitrary source point $z_1 \in \real^2 \setminus D$, with $E_n$ the element~(\ref{E_n}) from the fundamental sequence.  Here, $[\cdot]_1$ denotes the first column (the first sub-index of the boundary functions corresponds to the boundary curve where it is defined). The field
\begin{equation}\label{eq_exact}
u^{ex}_n (x) = [E_n (x,z_1)]_1, \quad x\in D,
\end{equation}
is then clearly the solution of \eref{s_e2}--\eref{s_bc3}. We consider as an approximate solution a solution of the form \eref{s_p}, where the densities $q^\ell_n$ for $\ell = 1,2$, and $n = 0,\ldots,N$, satisfy the system of equations
\begin{equation}\label{sie2}
\left\{\begin{array}{lr}
       \displaystyle{ -\frac{1}{2}q_n^1 (x)+\frac{1}{2\pi}\sum_{\ell=1}^2\int_{\Gamma_\ell}  T_xE_0(x, y)q_n^\ell(y)\, ds(y)  = G_{1,n} (x)},  & x\in\Gamma_1,\\
     \displaystyle{   \frac{1}{2}q_n^2(x)+\frac{1}{2\pi}\sum_{\ell=1}^2\int_{\Gamma_\ell} T_xE_0(x, y)q_n^\ell(y) \, ds(y) = G_{2,n} (x)}, & x\in\Gamma_2 ,
        \end{array}\right.
\end{equation}
with right-hand sides given here by
\begin{align*}
G_{1,n}(x) &= g_{1,n}(x) +
\frac{1}{2}\sum_{m=0}^{n-1}q_m^1(x) -\frac{1}{2\pi}\sum_{\ell=1}^2\sum_{m=0}^{n-1}\int_{\Gamma_\ell} T_x E_{n-m}(x,y)q_m^\ell(y)\, ds(y), \\
  G_{2,n}(x) &= g_{2,n}(x)-\frac{1}{2}\sum_{m=0}^{n-1}q_m^2(x)
 -\frac{1}{2\pi}\sum_{\ell=1}^2\sum_{m=0}^{n-1}\int_{\Gamma_\ell} T_x E_{n-m}(x,y)q_m^\ell(y)\, ds(y).
\end{align*} 

\begin{figure}[!ht]
\begin{center}
\includegraphics[scale=0.7]{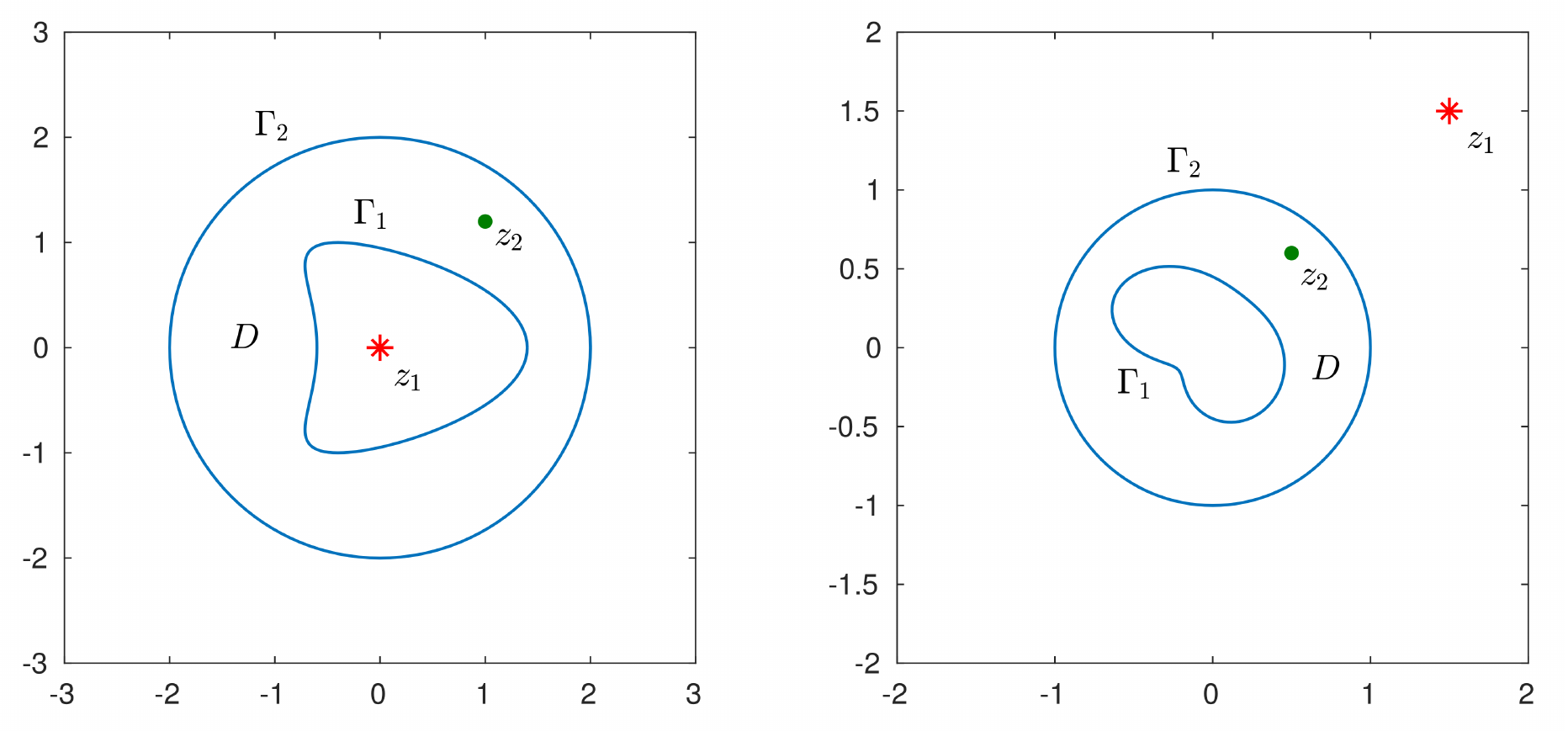}
\caption{The boundary $\Gamma_1 \cup \Gamma_2$ of $D$, the source point $z_1 \in \real^2 \setminus D$, and the measurement point $z_2 \in D$, in the first (left) and in the second (right) example.   }\label{fig1}
\end{center}
\end{figure}

The parametrized form of the approximate solution, considering \eref{s_p} and using the trapezoidal rule, is given by
\begin{equation}\label{eq_com}
u_n (x;M) = \frac{1}{2M} \sum_{\ell=1}^2 \sum_{m = 0}^n \sum_{k=0}^{2M-1} E _{n - m} (x, x_\ell (s_k)) \psi^\ell_m  (s_k),
\quad
x \in D,
\end{equation}
for $\psi_n^\ell(s)=q_n^\ell(x_\ell(s))|x'_\ell(s)|$. The numbers $\psi_n^\ell(s_k)$ are found from the corresponding parametric form of~(\ref{sie2}).

\noindent
{\bf Ex. 1}: In this first example, the outer boundary $\Gamma_2$ is set to be a circle with centre $(0,0)$ and radius 2, and the inner boundary $\Gamma_1$ is kite-shaped with parametrization
\[
\Gamma_1 = \{ x_1 (s) = (\cos s + 0.4 \cos 2s,\, \sin s): s\in [0,2\pi] \}.
\]
We consider the source point $z_1 = (0,\, 0)$ and the measurement point $z_2 = (1,\,1.2)$,  see the left picture in Fig. \ref{fig1}. The Lam\'e constants are $(\lambda, \, \mu) = (3,\, 2)$, and the density $\rho =1$. In Table \ref{table1}, we present the exact solution \eref{eq_exact} and the computed solution \eref{eq_com}, for $n=0,1,2$ and increasing number of mesh points $M$. We note that we get good accuracy with rather few points $M$. The order of error stays the same for $n=0,1,2$ and do not seem to grow much with increasing $n$ (remember that we have recurrence right-hand sides in~\eref{s_e2}--(\ref{s_bc3}) thus it is interesting to see how the error propagates with $n$). There is nothing special with the measurement point $z_2$, it can be moved around in the domain $D$ and even be taken on the boundary without affecting the results much. The Lam\'e constants can also be changed as can the source point; this will be done in the next example. 

\setlength{\tabcolsep}{18pt} 
\begin{table}[t]
\begin{center}
 \begin{tabular}{| c  | c  | c  | c  | } 
 \hline
 $M$ & $ (u_0)_1(z_2;M) $ & $(u_1)_1(z_2;M) $ & $(u_2)_1(z_2;M) $  
\\ \hline 
8 &    0.140558493353  & $   -0.054675857260  $ & $   -0.079614212909    $    \\ 
16 &  0.136678719983  & $   -0.053091680852   $ & $  -0.080635439439   $ \\
32 &    0.136669523731  & $   -0.053093812946  $ & $   -0.080639618965   $  \\
64 &   0.136669523108  & $   -0.053093813488   $ & $  -0.080639619432  $  \\
 \cline{1-1}\hhline{~===}

 \multicolumn{1}{c|}{} &  $ (u_0^{ex})_1 (z_2) $   & $ (u_1^{ex})_1 (z_2) $ & $ (u_2^{ex})_1 (z_2 ) $ \\
 \cline{2-4} 
\multicolumn{1}{c|}{} &       0.136669523108   & $  -0.053093813488    $ & $ -0.080639619432    $
 \\ \cline{2-4}
\end{tabular}
\vspace{0.4cm}
\caption{The first component of the computed (\ref{eq_com}) and the exact solution (\ref{eq_exact}) of \eref{s_e2}--\eref{s_bc3}, for the source point $z_1 = (0,\, 0)$ and the setup of the first example.}\label{table1}
\end{center}
\end{table}

\begin{figure}[!ht]
\begin{center}
\includegraphics[scale=0.75]{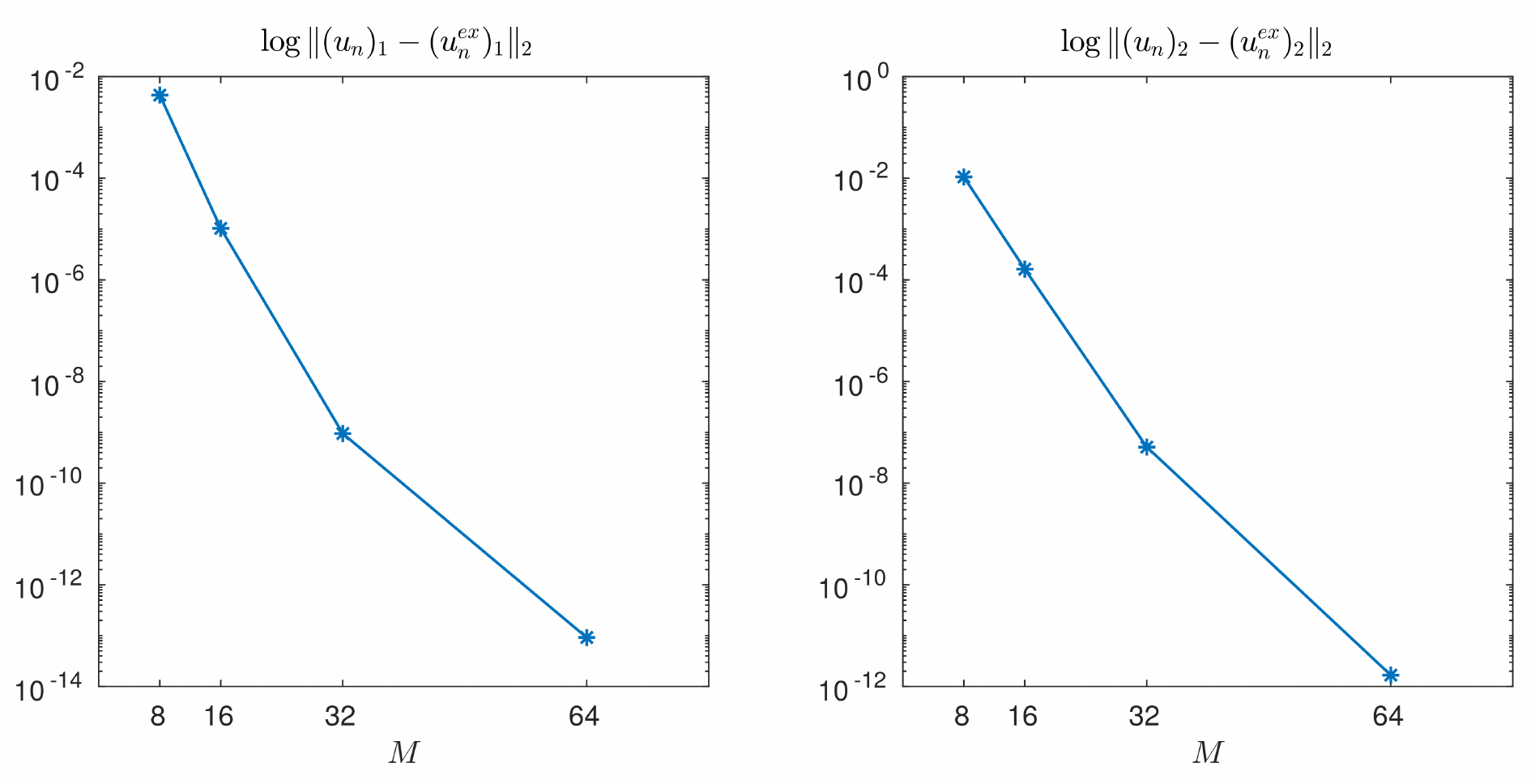}
\caption{The $L^2$-norm (logarithmic scale) of the difference between the computed (\ref{eq_com}) and the exact solution (\ref{eq_exact}). The convergence of the first component from Table \ref{table1} (left) and of the second component from Table \ref{table2} (right).  }\label{fig2}
\end{center}
\end{figure}

\noindent
{\bf Ex. 2}: We set the Lam\'e constants as $(\lambda, \, \mu) = (2,\, 1)$, and the density $\rho =1$. Both boundaries admit the form
\[
\Gamma_\ell = \{ x_\ell (s) = r_\ell (s)(\cos s,\, \sin s): s\in [0,2\pi] \}.
\]
We consider the radial functions
\[
r_1 (s) = \frac{0.9 + 0.6 \cos s -0.2 \sin 2s }{2+1.4 \cos s} \qquad \textrm{ and } \qquad  r_2 (s) = 1.
\]
The source point is now $z_1 = (1.5,\, 1.5)$ and we compute the fields at the measurement point $z_2 = (0.5,\,0.6)$, see the right picture in Fig. \ref{fig1}. Table \ref{table2} shows the exact solution \eref{eq_exact} and the computed solution \eref{eq_com}, for $n=0,1,2$ and varying number of mesh points $M$.  The
exponential convergence with respect to the spatial discretization is clearly
exhibited. This can also be seen in Fig. \ref{fig2}, where we plot the logarithm of the $L^2$-norm of the difference between the exact and the computed fields. This further justifies the claim in the previous example that the measurement point $z_2$ can be moved around within the solution domain and the similar accuracy will be obtained. Moreover, it shows also that the accuracy is similar for both components of the solution.  

\setlength{\tabcolsep}{18pt} 
\begin{table}[!ht]
\begin{center}
 \begin{tabular}{| c  | c  | c  | c  | } 
 \hline
 $M$ & $ (u_0)_2(z_2;M) $ & $(u_1)_2(z_2;M) $ & $(u_2)_2(z_2;M) $  
\\ \hline 
8 &    0.108260139657   &  0.026983579623  &   $-0.018597832984  $          \\ 
16 &   0.109123122407   &  0.030493156683  &   $-0.028643828117  $        \\
32 &    0.109244056837  &   0.030555054822  &  $ -0.028554171419 $  \\
64 &   0.109244013821   &  0.030555028775    & $-0.028554181488   $     \\
 \cline{1-1}\hhline{~===}

 \multicolumn{1}{c|}{} &  $ (u_0^{ex})_2 (z_2) $   & $ (u_1^{ex})_2 (z_2) $ & $ (u_2^{ex})_2 (z_2 ) $ \\
 \cline{2-4} 
\multicolumn{1}{c|}{} &       0.109244013821  &   0.030555028775    & $-0.028554181490     $
 \\ \cline{2-4}
\end{tabular}
\vspace{0.4cm}
\caption{The second component of the computed (\ref{eq_com}) and the exact solution (\ref{eq_exact}) of \eref{s_e2}--\eref{s_bc3}, for the source point $z_1 = (1.5,\, 1.5)$ and the setup of the second example.}\label{table2}
\end{center}
\end{table}

The obtained accuracy in the first two examples for direct problems serves as a numerical verification of the correctness of the derived formulas notably~\eref{TE_form}. We then turn to the ill-posed Cauchy problem~\eref{h_eq}--\eref{h_bc}.

\noindent
{\bf Ex. 3}: We first examine the feasibility of solving the sequence of ill-posed stationary problems \eref{s_e}--\eref{s_bc}. We consider again an arbitrary source point $z_1 \in \real^2 \setminus D$, and construct the boundary functions
\begin{equation}\label{f1e}
f_{2,n} (x)= [E_n (x,z_1)]_1 , \quad  g_{2,n}(x) = [T E_n (x,z_1)]_1, \quad x\in \Gamma_2 .
\end{equation}
We know then the exact solution \eref{eq_exact} and the computed \eref{eq_com}, where the densities now satisfies the linear system~\eref{sle}. Thus, we can compare the exact and the numerically calculated Cauchy data on the inner boundary $\Gamma_1$. We consider the boundary curves with  parametrization from the first example and the constants in the system as in the second example. The ill-posedness of the system is handled using Tikhonov regularization. The regularization parameter is chosen by trial and error. 
The source and the measurement points are $z_1 = (3,\,3)$ and $z_2 = \tfrac{\sqrt{2}}{2}(1,\,1) \in \Gamma_1 $, respectively.  The reconstructed Cauchy data are presented in Tables \ref{table3a} and \ref{table3b} for $n=0,5,10$. Here, we consider exact data and the regularization parameter is set to $10^{-10}$ for
$\kappa = 0.5$, and to $10^{-8}$ for $\kappa = 1$. As expected, the accuracy is slightly lower for the traction function since it involves derivatives of the solution. The scaling parameter $\kappa$ do influence the reconstructions but only slightly. Note that there is no dramatic increase in error as $n$ increases (recall that we have a recurrence right-hand side in~(\ref{sie}))

\setlength{\tabcolsep}{13pt} 
\begin{table}[!ht] 
\begin{center}
 \begin{tabular}{|c | c  | c  | c  | c  | } 
 \hline
$\kappa$ & $M$ & $ f_{1,0} (z_2) $ & $f_{1,5} (z_2) $ & $f_{1,10} (z_2) $
\\ \hline 
 \multirow{2}{*}{$0.5$} & 16 &  0.160838025793   &  0.028252624954     & 0.008865680398       \\
& 32 &     0.160981797423   &  0.028078500923  &   0.008781179908   \\ 
 \cline{1-2}\hhline{~~===}
 \multicolumn{2}{c|}{} &  $ [E_0 (z_2 ,z_1)]_1 $   & $ [E_5 (z_2 ,z_1)]_1 $ & $ [E_{10} (z_2 ,z_1)]_1 $ \\
 \cline{3-5} 
\multicolumn{2}{c|}{} & 0.160981796003  &   0.028078500985  &   0.008781181106 
 \\ \hhline{~~===} \cline{1-2}
 \multirow{2}{*}{1} & 16 &   0.035987791353  &   0.008580531214    & $-0.037640150526 $  \\
& 32 &   0.036002107356  &   0.008720257700   &  $\phantom{-}0.002279598771  $   \\
 \cline{1-2}\hhline{~~===}
 \multicolumn{2}{c|}{} &  $ [E_0 (z_2 ,z_1)]_1 $   & $ [E_5 (z_2 ,z_1)]_1 $ & $ [E_{10} (z_2 ,z_1)]_1 $ \\
 \cline{3-5} 
 \multicolumn{2}{c|}{} & 0.036002151515   &  0.008720380239    & $\phantom{-}0.002279504879 $
 \\ \cline{3-5}
\end{tabular}
\caption{The reconstructed (via~(\ref{f1})) and the exact~(\ref{f1e}) values of the boundary function $f_{1,n}$ on $\Gamma_1$, for the source point $z_1 = (3,\,3)$, and the setup of the third example.}\label{table3a}
\end{center}
\end{table}

\begin{table}[!ht] 
\begin{center}
 \begin{tabular}{|c | c  | c  | c  | c  | } 
 \hline
$\kappa$ & $M$ & $ g_{1,0} (z_2) $ & $g_{1,5} (z_2) $ & $g_{1,10} (z_2) $
\\ \hline 
 \multirow{2}{*}{$0.5$} & 16 &  $ -0.158494501053  $&   0.088261735896  &$    -0.031767212741  $    \\
& 32 &   $   -0.156020048879 $&    0.084724640114  &$   -0.031135310036  $   \\ 
 \cline{1-2}\hhline{~~===}
 \multicolumn{2}{c|}{} &  $ [TE_0 (z_2 ,z_1)]_1 $   & $ [TE_5 (z_2 ,z_1)]_1 $ & $ [TE_{10} (z_2 ,z_1)]_1 $ \\
 \cline{3-5} \multicolumn{2}{c|}{} & $-0.156020119157  $&    0.084724639664  &$    -0.031135305722  $
 \\ \hhline{~~===} \cline{1-2}
 \multirow{2}{*}{1} & 16 &  $ -0.069488888403  $ & $   -0.022080526430    $ & $ -0.767369563232 $   \\
& 32 &  $ -0.069308164504 $ &  $ -0.017890055807 $ & $  -0.016726813312$    \\
 \cline{1-2}\hhline{~~===}
 \multicolumn{2}{c|}{} &  $ [TE_0 (z_2 ,z_1)]_1 $   & $ [TE_5 (z_2 ,z_1)]_1 $ & $ [TE_{10} (z_2 ,z_1)]_1 $ \\
 \cline{3-5} 
 \multicolumn{2}{c|}{} & $-0.069309067338$ & $      -0.017888657343 $ & $      -0.016727744904  $
 \\ \cline{3-5}
\end{tabular}
\caption{The reconstructed (via~(\ref{g1})) and the exact values~(\ref{f1e}) of the boundary function $g_{1,n}$ on $\Gamma_1$, for the source point $z_1 = (3,\,3)$, and the setup of the third example.}\label{table3b}
\end{center}
\end{table}

To investigate the stability of the method against noise, we add noise to the boundary function $g_{2,n} \textrm{ on } \Gamma_2$, with respect to the $L^2-$norm 
\[
g^\delta_{2,n}  = g_{2,n} + \delta \frac{\|g_{2,n}\|_2}{\|v \|_2}v,
\]
for a given noise level $\delta$ and a normally distributed random variable $v$. We define the following relative $L^2$-errors on the inner boundary 
\begin{equation}\label{e2f}
e^2_f (n) = \frac{\displaystyle\int_0^{2\pi}\left( f_{1,n} (x_1 (s)) - [E_n (x_1 (s),z_1)]_1 \right)^2\hspace*{0.05cm} ds}{\displaystyle\int_0^{2\pi} ([E_n (x_1 (s),z_1)]_1 )^2\, ds} ,
\end{equation}
and
\begin{equation}\label{e2g}
e^2_g (n) = \frac{\displaystyle\int_0^{2\pi}\left( g_{1,n} (x_1 (s)) - [TE_n (x_1 (s),z_1)]_1 \right)^2\hspace*{0.05cm} ds}{\displaystyle\int_0^{2\pi} ([TE_n (x_1 (s),z_1)]_1 )^2\, ds} .
\end{equation}
The error terms~(\ref{e2f})--(\ref{e2g}) are numerically calculated using the trapezoidal rule. In Table \ref{table4} are 
the values of the error terms for exact ($\delta=0\%$) and noisy data ($\delta=3\%$) with regularization parameter $10^{-8}$ and $10^{-2}$, respectively. We have set $\kappa=1$ and $M=32$ and we keep the same source point $z_1$.
We observe that the errors do not increase dramatically when $n$ increases.  Note that keeping increasing $M$ will not decrease errors much further due to the increasing ill-conditioning of the linear system, which in turn is due to the ill-posedness of the underlying elastodynamic problem.

\setlength{\tabcolsep}{14pt} 
\begin{table}[t]
\begin{center}
 \begin{tabular}{| c  | c  | c  | c  |  c |} 
 \cline{2-5}
  \multicolumn{1}{c|}{} &  \multicolumn{2}{c|}{$\delta = 0\%$} &  \multicolumn{2}{c|}{$\delta = 3\%$}
\\ \hline 
$n$ & $e_f (n)$ & $e_g (n)$ & $e_f (n)$ & $e_g (n)$  \\ \hline
5 & 5.9908E$-$06  & 7.7404E$-$05    &     0.14558   &   0.71139   \\ \hline
10 &    1.0238E$-$05  & 8.0608E$-$05  &    0.14038   &   0.65011\\ \hline
15 & 1.0642E$-$05 &   4.6722E$-$05 &      0.39071    &   0.64990\\ \hline
20 & 7.8734E$-$05 & 5.5732E$-$04 &   0.41676  & 0.82314\\ \hline
\end{tabular}
\vspace{0.4cm}
\caption{The errors~(\ref{e2f})--(\ref{e2g}) in the third example with exact ($\delta=0\%$) and noisy data ($\delta=3\%$) for the source point $z_1 = (3,\,3)$, and $M=32$.}\label{table4}
\end{center}
\end{table}

Considering the expansion \eref{lag6}, we can also compare the time-dependent computed solution
\begin{equation}\label{un}
\tilde u (x,t) = \kappa \sum_{n=0}^{N-1} u_n (x;M) L_n (\kappa t), \quad x\in D,
\end{equation}
with the fundamental solution (truncated form)
\begin{equation}\label{Ext}
\tilde E (x, t) = \kappa \sum_{n=0}^{N-1} [E_n (x,z_1)]_1 L_n (\kappa t), \quad x\in D.
\end{equation}
The values of the first component of the computed solution are presented in Table \ref{table5} for varying $N$ and $M$ and at different times $t$. The values at the interior point $z_2 = (1,\,1.2)$ for $\kappa=1$, are given for different regularization parameters depending on the error level and the number of Fourier coefficients. For noise free data, we use $10^{-6}$ as regularization parameter and for $\delta=3\%$ we set it to $10^{-4}$.

\setlength{\tabcolsep}{6pt} 
\begin{table}[!ht]
\begin{center}
 \begin{tabular}{| c | c  | c  | c  |  c | c|}  
  \cline{3-6}
  \multicolumn{2}{c|}{} &  \multicolumn{2}{c|}{$N=15$} &  \multicolumn{2}{c|}{$N=20$}
\\ \hline 
$t$ & $M$ &  $\delta = 0\%$& $\delta = 3\%$  & $\delta = 0\%$ & $\delta = 3\%$
\\ \hline 
 \multirow{3}{*}{1} & 16 &   0.46387720765  &  0.56160816191  & 0.34605948854 &  0.54725934363 
 \\
 & 32 &  0.54302224375    &  0.54522283409 &  0.54096499271  &  0.53709220126  \\ 
 & $[\tilde E (t)]_1$ & \multicolumn{2}{c|}{\cellcolor[gray]{0.8}   0.54305542205}  &  \multicolumn{2}{c|}{\cellcolor[gray]{0.8}  0.54099818785}   \\
 \cline{1-1}\hhline{======}
  \multirow{3}{*}{2} & 16 &   0.59385442392    &  0.46302021628 &  $-1.41326137373$ &  0.36623652546    \\
 & 32 &   0.48564171165   & 0.48223841917   &   $\phantom{-}0.47964119241$&  0.47575211502  
 \\ 
 & $[\tilde E (t)]_1$ & \multicolumn{2}{c|}{\cellcolor[gray]{0.8} 0.48567189228}     & \multicolumn{2}{c|}{\cellcolor[gray]{0.8}  0.47968395291   } \\
 \cline{1-1}\hhline{======}
  \multirow{3}{*}{3} & 16 & 0.02702301633   &  0.18548879529 &  $-0.35158612858 $ &  0.18173428987     \\
 & 32 & 0.16541236649    & 0.17064119470 &  $\phantom{-}0.15815935533$ & 0.15671723036     \\ 
 & $[\tilde E (t)]_1$ & \multicolumn{2}{c|}{\cellcolor[gray]{0.8}    0.16543897091}  & \multicolumn{2}{c|}{\cellcolor[gray]{0.8}  0.15818608537 } \\ \hline
\end{tabular}
\vspace{0.4cm}
\caption{Numerical values of the first component of the computed solution $\tilde{u} (z_2,t)$ from~(\ref{un})  and the exact solution $\tilde E (z_2,t)$ (rows in grey) from~(\ref{Ext}) in the third example, at the measurement point $z_2 =  (1,\, 1.2)$, for exact ($\delta=0\%$) and noisy ($\delta=3\%$) data and various times $t$.}\label{table5}
\end{center}
\end{table}

\noindent\\
{\bf Ex. 4}: In this example, we examine solving (\ref{h_eq})--(\ref{h_bc}) using our method to reconstruct in a stable way the missing lateral Cauchy data on $\Gamma_1$, but we have no exact solution to test against. Instead, we numerically construct the Cauchy data on $\Gamma_2$ by solving the well-posed problem
\begin{subequations}
\begin{alignat}{2}
\frac{{\partial^2 u}}{{\partial \,t^2}} &= \Delta^* u,\quad & \textrm{in }\,\, & D \times ({0,\infty} ),  \label{third2}\\
u&= f_1, & \textrm{on }\,\, & \Gamma_1 \times (0,\infty), \label{third_bc2}\\
u&= f_2, & \textrm{on }\,\, & \Gamma_2 \times (0,\infty), \label{third_bc3}
\end{alignat}
\end{subequations}
together with homogeneous initial conditions. We consider the boundary functions
\[
f_1 (x,t) = f(t) (1,\,1)^\top , \quad f_2 (x,t) = (0,\,0)^\top , \quad \mbox{with} \quad f(t) = \frac{t^2}{4}e^{-t+2},
\]
which admits the expansion
\begin{equation}\label{expansion}
f(t) = \frac{\kappa e}{4} \sum_{n=0}^\infty \frac{2+\kappa n (\kappa (n-1)-4)}{(\kappa+1)^{n+3}} L_n (\kappa t).
\end{equation}

We solve the direct problem for $(\lambda,\mu) = (2,\,1), \, \rho = 1$ and $\kappa =1$, using the boundary curves with parametrization as in the first example and $M=64$. In order to avoid an ``inverse crime'', we solve the inverse problem using half the number of collocation points and adding $3\%$ noise on the computed function $g_2$. The exact time dependent function $f_1(x,t)$ is constructed using the above expansion truncated at $N=20$. Table \ref{table6} shows the first component of the reconstructed functions on $\Gamma_1$ for the sequence of stationary problems~\eref{s_e}--\eref{s_bc} and the corresponding exact value (obtained by the Laguerre transform of $f_1$). Table \ref{table7} shows the reconstructed time dependent function on $\Gamma_1$, which we compare with the above exact value $f_1$ (expansion truncated at $N=20$). The measurement point is $z_2 = (1.4,\,0)$ and for $N=10,\,15$ and $20$, we set the regularization parameter to $10^{-4}$ for the first two cases and $10^{-3}$ for the latter. 

We define the transient error term
\begin{equation*}
e^2  = \frac{\displaystyle\int_0^T\int_0^{2\pi}\left( f_1 (x_1 (s),t) - \tilde u (x_1 (s),t) \right)^2\, ds\hspace*{0.05cm}dt}{\displaystyle\int_0^T\int_0^{2\pi} f_1 (x_1 (s),t)^2\, ds\hspace*{0.05cm}dt} ,
\end{equation*}
where $\tilde u$ represents the solution of the inverse problem restricted to the inner boundary. The errors with respect to $N$ for $T=3$ are shown in the last column of Table \ref{table7}.  For the time integration we use the trapezoidal rule with step size $0.2$.

\begin{table}[!ht] 
\begin{center}
 \begin{tabular}{|c | c  | c  | c  | c  | } 
 \hline
$N$ & $f_{1,0} (z_2)$ & $ f_{1,5} (z_2) $ & $g_{1,0} (z_2) $ & $g_{1,5} (z_2) $
\\ \hline 
10 &   0.162104216452 &    $-0.012768621007$ &   1.129505738211    & $-0.069348086158 $ \\ \hline
15 & 0.167765321801  &   $-0.011047131720$  &  1.186907689072    & $-0.053265077003$ \\ \hline
20 & 0.171214630252  &   $-0.010092156278 $ & 1.194786602851    & $-0.042902486556$  \\ \cline{1-1}\hhline{~====} 
\multicolumn{1}{c|}{} & \cellcolor[gray]{0.8}   0.169892614279    & \cellcolor[gray]{0.8} $ -0.010618288392 $ &  \cellcolor[gray]{0.8} 1.208462595905  &  \cellcolor[gray]{0.8}  $-0.044010338751$    \\ \cline{2-5}
\end{tabular}
\caption{The reconstructed time independent boundary functions (generated via~(\ref{f1})--(\ref{g1})) for noisy data ($\delta=3\%$) at the measurement point $z_2 = (1.4, \,0)$ and $M=64$. The row in grey shows the exact values given from the direct problem \eref{third2}--\eref{third_bc3}. }\label{table6}
\bigskip
 \begin{tabular}{|c | c  | c  | c  |  c |} 
 \hline
$N$ & $f_1 (z_2 ,1)$ & $f_1 (z_2 ,2)$ & $f_1 (z_2 ,3)$  & $e$ 
\\ \hline 
10 & 0.236547147854275 &  0.361215644486791 &  0.280287929224566 & 0.12845 \\ \hline
15 &    0.249079092594895  & 0.360552859388251 &  0.320751004460423 & 0.14661\\ \hline
20 & 0.249089211842107 &  0.368579099163008 &  0.317125178566481 & 0.30520
\\\cline{1-2}\cline{4-5}
\hhline{~===~} 
\multicolumn{1}{c |}{} & \cellcolor[gray]{0.8}      0.250013290485691  & \cellcolor[gray]{0.8}  0.367850715255786  & \cellcolor[gray]{0.8}  0.304541324965188    \\ \cline{2-4}
\end{tabular}
\caption{The reconstructed time dependent boundary function (via (\ref{lag6}) and (\ref{eq_com})) on $\Gamma_1\times (0,T)$, with $T=3$, for noisy data ($\delta=3\%$) at the measurement point $z_2 = (1.4, \,0)$, for $t=1,2,3$ and $M=64$. The row in grey shows the exact values given from the truncated expansion of $f(t)$ in~(\ref{expansion}) calculated for $N=20$. The transient error is shown in the last column. }\label{table7}
\end{center}
\end{table}

\section{Conclusion} 
A  boundary integral equation approach has been developed for an ill-posed later Cauchy problem for the elastodynamic equation. Given the function value and traction on the outer boundary of an annular domain the corresponding data are reconstructed on the inner domain. Applying the Laguerre transformation in time the problem is reduced to a sequence of ill-posed stationary problems. Using what is known as a fundamental sequence to these stationary equations, the solutions are represented as  single-layer potentials over the boundary without involving any domain integrals. Moreover, exact representations in terms of boundary integrals were given for the missing data on the inner boundary. Careful analysis of the singularities of the kernels of the boundary integrals makes it possible to introduce a suitable splitting such that efficient quadrature methods can be applied for numerical discretisations. Numerical experiments included show that accurate solutions can be obtained to the time-dependent system both for direct and ill-posed problems.


\end{document}